\documentclass[11pt]{article}
\usepackage{amsfonts}
\usepackage{amsmath}
\usepackage{amsthm}
\usepackage{url}
\usepackage{amssymb}
\usepackage{latexsym}
\usepackage[square, comma, sort&compress]{natbib}
\usepackage{float}
\usepackage{labelfig}
\usepackage{epsfig,psfig}
\usepackage{epstopdf}

\usepackage{graphics}

\def\pg{\mathcal{P}}
\def\S{\Sigma}

\def\kk{\kappa}
\def\gg{\gamma}
\def\aa{\alpha}
\def\bb{\beta}
\def\pp{\phi}

\def\Mod{{\rm{Mod}}}
\def\st{{\rm{St}}}
\long\def\symbolfootnote[#1]#2{\begingroup\def\thefootnote{\fnsymbol{footnote}}
\footnote[#1]{#2}\endgroup}

\newtheorem*{thma}{Theorem A}
\newtheorem*{thmb}{Theorem B}
\newtheorem*{thmc}{Theorem C}

\newtheorem{conjecture}{Conjecture}

\newtheorem{corollary}[conjecture]{Corollary}
\newtheorem{definition}[conjecture]{Definition}

\newtheorem{lemma}[conjecture]{Lemma}

\newtheorem{theorem}[conjecture]{Theorem}

\begin{document}

\begin{center}

\bf{SIMPLICIAL EMBEDDINGS BETWEEN PANTS GRAPHS} \\

\bigskip

Javier Aramayona \footnote{Supported by M.E.C. grant MTM2006/14688 and NUI Galway's Millennium Research Fund.} \\

\end{center}


\begin{abstract}

\noindent We prove that, except in some low-complexity cases, every locally injective simplicial map between pants graphs is induced by a $\pi_1$-injective embedding between the corresponding surfaces. 

\end{abstract}


\section{Introduction and main results}
\label{intro}

\symbolfootnote[0]{Date: \today}

To a surface $\S$ one may associate a number of naturally defined objects -- its Teichm\"uller space, mapping class group, curve or pants graph, etc. An obvious problem is then to study embeddings between objects in the same category, where the term ``embedding" is to be interpreted suitably in each case, for instance ``isometric embedding" in the case of Teichm\"uller spaces, ``injective homomorphism" in the case of mapping class groups, and ``injective simplicial map" in the case of curve and pants graphs. 

For pants graphs, this problem was first studied by D. Margalit \cite{margalit}, who showed that every
automorphism of the pants graph is induced by a self-homeomorphism of $\S$. More concretely, let ${\rm{Mod}}(\S)$ be the mapping class group of $\S$, which acts on the pants graph $\pg(\S)$ by simplicial automorphisms, and let ${\rm{Aut}}(\pg(\S))$ be the group of all simplicial automorphisms of $\pg(\S)$. Let $\kk(\S)$ be the complexity of $\S$, that is, the cardinality of a pants decomposition of $\S$. The following is part of Theorem 1 of \cite{margalit}:

\begin{theorem}[\cite{margalit}] 
If $\S$ is a compact, connected, orientable surface with $\kk(\S)>0$, then the natural homomorphism ${\rm{Mod}}(\S) \rightarrow {\rm{Aut}}(\pg(\S))$ is surjective. Moreover, if $\kk(\S) > 3$ then it is an isomorphism.
\label{dan}
\end{theorem}

The main purpose of this note is to extend Margalit's result to (locally) injective simplicial maps between pants graphs. We note that examples of such maps are plentiful. Indeed, let $\S_1$ be an essential subsurface of $\S_2$ (see Section \ref{defs} for definitions) whose every connected component  has positive complexity. Then
one may construct  an injective simplicial map $\pp: \pg(\S_1) \rightarrow \pg(\S_2)$ by first choosing a multicurve $Q$ that extends any pants decomposition of $\S_1$ to a pants decomposition of $\S_2$ and then setting $\phi(v) = v \cup Q$. 

Our main result asserts that, except in some low-complexity cases, this is the only way in which injective simplicial maps of pants graphs arise. Given a simplicial map $\pp: \pg(\S_1) \rightarrow \pg(\S_2)$ and a $\pi_1$-injective embedding $h: \S_1 \rightarrow \S_2$,  we say that $\pp$ is {\em{induced}} by $h$ if there exists a multicurve $Q$ on $\S_2$, disjoint from $h(\S_1)$, such that $\pp(v) = h(v) \cup Q$ for all vertices $v$ of $\pg(\S_1)$. In particular, $Q$ has cardinality $\kk(\S_2) - \kk(\S_1)$.  We will show:

\begin{thma}
Let $\S_1$ and $\S_2$ be compact orientable surfaces of negative Euler characteristic, such that each connected component of $\S_1$ has complexity at least 2. Let $\pp: \pg(\S_1) \rightarrow \pg(\S_2)$ be an injective simplicial map. Then there exists a $\pi_1$-injective embedding $h: \S_1 \rightarrow \S_2$  that induces $\pp$.

\label{main-theorem}
\end{thma}

We note that the hypothesis that all connected components of $\S_1$ have complexity at least 2 is necessary, since the pants graph of the 1-holed torus and the 4-holed sphere are isomorphic (see \cite{minsky}, for instance). 

\bigskip

\noindent{\bf{Remark.}}  In the case of  curve graphs, Teichm\"uller spaces and mapping class groups, there exist embeddings for which there are no $\pi_1$-injective embeddings of the corresponding surfaces. First, one may construct an injective simplicial map from the curve graph of a closed surface $X$ to that of $X-p$, by considering a point $p$ in the complement of the union of all simple closed geodesics on 
$X$. Next, any finite-degree cover $\tilde{Y} \to Y$  gives rise to an isometric embedding $T(Y) \to T(\tilde{Y})$ of Teichm\"uller spaces, so we may take $Y$ to be a closed surface in order to produce the desired example. Finally, there exist injective homomorphisms of mapping class groups with no $\pi_1$-injective embeddings between the corresponding surfaces, se \cite{birman-hilden} and \cite{ALS}.

\bigskip

 In order to prove Theorem A, we will closely follow Margalit's strategy in \cite{margalit} for proving Theorem 1. In Section \ref{defs} we will introduce the pants graph and its natural subgraphs. In Section \ref{objects} we will study some objects in the pants graph, namely {\em{Farey graphs}} and {\em{admissible tuples}}, which appear, or at least have their origin, in \cite{margalit}. Most importantly, the structure of these objects is preserved by injective simplicial maps. Using this, in Section \ref{proof}  we will show the following result, which will constitute the main step for proving Theorem A.  Given a multicurve $Q$ on $\S$, let $\pg_Q$ be the subgraph of $\pg(\S)$ spanned by those vertices containing $Q$. A {\em{non-trivial component}} of  $\S$ is a connected component of $\S$ of positive complexity.

\begin{thmb}
Let $\S_1$ and $\S_2$ be compact orientable surfaces such that every connected component of $\S_1$ has positive complexity.  Let $\pp: \pg(\S_1) \rightarrow \pg(\S_2)$ be an injective simplicial map. Then the following hold:

\begin{enumerate}

\item[(1)] $\kk(\S_1) \leq \kk(\S_2)$,

\item[(2)] There exists a multicurve $Q$ on $\S_2$, of cardinality $\kk(\S_2) - \kk(\S_1)$, such that $\pp(\pg(\S_1))= \pg_Q$. In particular, $\pg(\S_1) \cong \pg(\S_2 -Q)$;

\item[(3)] $\S_1$ and $\S_2 -Q$ have the same number of non-trivial components. Moreover, if $X_1, \ldots, X_r$ and $Y_1, \ldots, Y_r$ are, respectively, the non-trivial components of $\S_1$ and $\S_2 - Q$ then, up to reordering the indices, $\pp$ induces an isomorphism $\pp_i: \pg(X_i) \rightarrow \pg(Y_i)$. In particular, $\kk(X_i) = \kk(Y_i)$.

\end{enumerate}

\label{expanded}
\end{thmb}

Theorem B itself has an interesting consequence for pants graph automorphisms. More concretely, in Corollary \ref{strata} of Section \ref{apps} we will see that pants graph automorphisms preserve the pants graph stratification (see Section \ref{defs} for definitions). This  implies Theorem \ref{dan} if  $\S$ is not the 2-holed torus. The case of the 2-holed torus needs some extra care but it also follows from Theorem B by applying the same strategy of \cite{margalit}, Section 5.

Finally, in Section \ref{final} we will prove Theorem A, which will follow easily from Theorem B and  the (folklore) classification of pants graphs up to isomorphism, included as Lemma \ref{isom-class} in Section \ref{final}. 

We remark that, even though Theorems A and B are stated for injective simplicial maps, our arguments will only require the maps to be simplicial and {\em{locally injective}}, that is, injective on the star of every vertex of $\pg(\S_1)$. The {\em{star}} of a vertex is defined as the union of all edges incident on it. In particular, we have:

\begin{thmc}
Let $\S_1$ and $\S_2$ be compact orientable surfaces of negative Euler characteristic, such that each connected component of $\S_1$ has complexity at least 2. Let $\pp: \pg(\S_1) \rightarrow \pg(\S_2)$ be a locally injective simplicial map. Then there exists a $\pi_1$-injective embedding $h: \S_1 \rightarrow \S_2$  that induces $\pp$.
\label{local-main-theorem}
\end{thmc}

Finally, we point out that a number of authors have studied embeddings in the context of mapping class groups and other complexes associated to surfaces. References include  \cite{ALS}, \cite{behrstock-margalit}, \cite{bell-margalit}, \cite{birman-hilden}, \cite{irmak1}, \cite{irmak2}, \cite{irmak3}, \cite{irmak-korkmaz}, \cite{irmak-mccarthy}, \cite{ivanov}, \cite{ivanov-mccarthy}, \cite{korkmaz}, \cite{luo}, \cite{paris-rolfsen}, \cite{schmutz}, \cite{shackleton}.

\bigskip

\noindent{\bf{Acknowledgements.}} Parts of this work were completed during visits of the author to 
MSRI and UIUC. I would like to express my gratitude towards these institutions. I warmly thank
Jeff Brock, Cyril Lecuire, Chris Leininger, Dan Margalit, Kasra Rafi, Ken Shackleton, Juan Souto and Ser Peow Tan for conversations, and  many invaluable comments and suggestions.


\section{Definitions and basic results}
\label{defs}

\subsection{Surfaces and curves.}

Let $\S$ be a compact orientable surface whose every connected component has negative Euler characteristic.  If $g$ and $b$ are, respectively, the genus and number of boundary components of $\S$, we will refer to the number $\kk(\S) = 3g -3+b$ as the {\em{complexity}} of $\S$. As mentioned in Section \ref{intro}, a {\em{non-trivial component}} of $\S$ is a connected component of $\S$ that has positive complexity. 

A subsurface $X \subset \S$ is said to be {\em{essential}} if no components of $X$ are parallel to $\partial \S$ and every component of  $\partial X$ determines either a non-homotopically trivial simple closed curve on $\S$ or a component of $\partial \S$. \break Throughout this note we will only consider essential subsurfaces whose every connected component has negative Euler characteristic. Two subsurfaces are said to be {\em{disjoint}} if they can be homotoped away from each other. 

A simple closed curve on $\S$ is said to be {\em{peripheral}} if it is homotopic to a component of $\partial \S$. By a {\em{curve}} on $\S$ we will mean a homotopy class of non-trivial  and non-peripheral simple closed curves on $\S$. The {\em{intersection number}} between two curves $\aa$ and $\bb$ is defined as $$i(\aa,\bb) = \text{min} \{ | a \cap b | : a \in \aa, b\in \bb\}.$$ 

If $i(\aa,\bb)=0$, we say that $\aa$ and $\bb$ are {\em{disjoint}}.  A {\em{multicurve}} is a collection of distinct and disjoint curves on $\S$.  Given a multicurve $Q$ on $\S$,  the {\em{deficiency}} of $Q$ is defined to be $\kk(\S) - |Q|$. A {\em{pants decomposition}} is a multicurve of cardinality $\kk(\S)$ (and so maximal with respect to inclusion). Note, if $Q$ is a pants decomposition then $\S -Q$ is a disjoint union of 3-holed spheres, or {\em{pairs of pants}}.

If $X \subset \S$ is a (not necessarily proper) subsurface, we say that a collection $\cal{A}$ of curves on $X$ {\em{fills}} $X$ if, for every curve $\gamma$ on $X$, there exists $\aa \in \cal{A}$ with $i(\aa, \gg)>0$. 
 In particular, if $\kk(X) = 1$ then any pair of distinct curves fill $X$.


\subsection{The pants graph.}

\begin{figure}[h]
\begin{center}
\AffixLabels{\centerline{\epsfig{file = 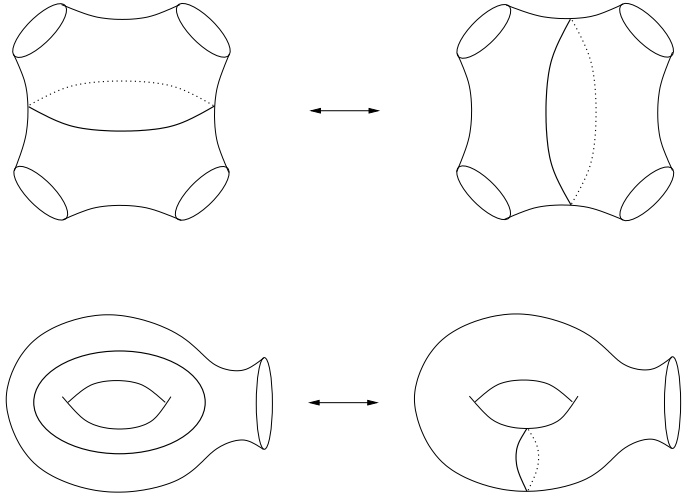, width=9cm, angle= 0}}}
\end{center}
\caption{The two types of elementary move.}
\label{fig:elementarymoves}
\end{figure}

 We say that two pants decompositions are related by an {\em{elementary move}} if they have a deficiency 1 multicurve in common, and the remaining
two curves either fill a 4-holed sphere and intersect exactly twice, or they fill a 1-holed torus and intersect exactly once. See  Figure \ref{fig:elementarymoves}.

The {\em{pants graph}} $\pg(\S)$ of $\S$ is the simplicial graph whose vertex set is the set of all pants decompositions of $\S$ and where two vertices are connected by an edge if the corresponding pants decompositions are related by an elementary move. A {\em{path}} in $\pg(\S)$ is a sequence $v_1, \ldots, v_n$ of adjacent vertices of $\pg(\S)$. A {\em{circuit}} is a path  $v_1, \ldots, v_n$ such that $v_1 = v_n$ and $v_i \neq v_j$ for all other $i,j$.

The pants graph was introduced by Hatcher-Thurston in \cite{hatcher-thurston}, who proved it is connected (see the remark on the last page of \cite{hatcher-thurston}). A detailed proof was then given by Hatcher-Lochak-Schneps in \cite{hls}, where they proved that attaching 2-cells to  finitely many types of circuits in $\pg(\S)$ produces a simply-connected 2-complex, known as the {\em{pants complex}}.
The graph $\pg(\S)$ becomes a geodesic metric space by declaring each edge to have length 1, and Brock \cite{brock} recently showed that $\pg(\S)$ is quasi-isometric to the Weil-Petersson metric on the Teichm\"uller space of $\S$.


\subsection{Natural subgraphs.} 

As mentioned in the introduction, if $Y \subset \S$ is an essential subsurface with no trivial components, then the inclusion map induces  an injective simplicial map $\pg(Y) \to \pg(\S)$, and so we can regard $\pg(Y)$ as a connected subgraph of $\pg(\S)$. If $Y_1, Y_2 \subset \S$ are disjoint essential subsurfaces of positive complexity, then $\pg(Y_1) \times \pg(Y_2)$, the 1-skeleton of the product of $\pg(Y_1)$ and $\pg(Y_2)$, is a connected subgraph of $\pg(\S)$. Moreover, if $\S$ is not connected, then $\pg(\S)$ is the 1-skeleton of the product of the pants graphs of its non-trivial components. 

Given a multicurve $Q$ on $\S$, let $\pg_Q$ be the subgraph of $\pg(\S)$ spanned by those vertices of $\pg(\S)$ that contain $Q$. It will be convenient  to consider the empty set as a multicurve, in which case we set $\pg_\emptyset$ to be equal to $\pg(\S)$. Note that $\pg_Q$ is connected for all multicurves $Q$; indeed, if $Q$ is strictly contained in a pants decomposition then $\pg_Q$ is naturally isomorphic to $\pg(\S -Q)$, and if $Q$ is itself a pants decomposition then $\pg_Q$ is equal to $Q$. 

If $Q_1$ and $Q_2$ are multicurves on $\S$, then $\pg_{Q_1}\cap \pg_{Q_2}\neq \emptyset$ if and only if  $Q_1 \cup Q_2$ is a multicurve, in which case $\pg_{Q_1} \cap \pg_{Q_2}  = \pg_{Q_1 \cup Q_2}$. Furthermore, $\pg_{Q_1} \subset \pg_{Q_2}$ if and only if $Q_2 \subset Q_1$.  This endows the pants graph with a {\em{stratified structure}},  analogous to the stratification of the Weil-Petersson completion (see \cite{wolpert}), with  strata all the subgraphs of the form $\pg_Q$, for some multicurve $Q$. Then $\pg(\S)$ is the union of all  strata, and two strata intersect over a stratum if at all.


\section{Some objects in the pants graph.}
\label{objects}

\subsection{Farey graphs.}

The {\em{standard Farey graph}} is the simplicial graph whose vertex set is $\mathbb{Q} \cup \{\infty\}$ and where two vertices $p/q$ and $r/s$, in lowest terms,
are connected by an edge if $|ps - rq | = 1$. It is usually represented as an ideal triangulation of the Poincar\'e disc model of  the hyperbolic plane. By a {\em{Farey graph}} we will mean an isomorphic copy of the standard Farey graph. The following result is implicit in the proof of Lemma 1 of \cite{margalit}.

\begin{lemma}[Structure of Farey graphs]
A subgraph $F$ of $\pg(\S)$ is a Farey graph if and only if $F = \pg_Q$, for some deficiency 1 multicurve $Q$.
\label{farey}
\end{lemma}

\noindent{\bf{Proof.}} Consider $\pg_Q$, where $Q$ is a deficiency 1 multicurve. Then $\S -Q$ has a unique non-trivial component $X$, which has complexity 1, and so it is either a 1-holed torus or a 4-holed sphere. In either case, $\pg(X)$ is a Farey graph (see, for instance, \cite{minsky}, Section 3), and $\pg_Q \cong \pg(\S - Q) \cong \pg(X)$.

Conversely, if $\Delta \subset \pg(\S)$ is a circuit of length 3, then $\Delta \subset \pg_Q$ for some deficiency 1 multicurve $Q$, by the definition of adjacency in $\pg(\S)$.  The result now follows easily from the observation that any two vertices of a Farey graph can be connected by a sequence of circuits of length 3 such that any two consecutive such circuits have exactly two vertices in common. $\square$

\bigskip

In the situation of Lemma \ref{farey}, we will say that $F$ is {\em{determined}} by $Q$. Let $e$ be an edge of $\pg(\S)$, and let $u$ and $v$ be its endpoints. Then $e$ is contained in a unique Farey graph, determined by the deficiency 1 multicurve $u \cap v$.  Given a vertex $u$ of $\pg(\S)$, observe that there are exactly $\kk(\S)$ distinct Farey graphs containing $u$, determined
 by the $\kk(\S)$ distinct deficiency 1 multicurves contained in $u$.  We state this observation as a separate lemma, as we will make extensive use of it later.

\begin{lemma}
Given any vertex $u$ of $\pg(\S)$, there are exactly $\kk(\S)$ distinct Farey graphs containing $u$. $\square$
\label{pinwheel}
\end{lemma}

As mentioned in Section \ref{intro}, the {\em{star}} $\st(u)$ of a vertex $u$ of $\pg(\S)$ is the union of all edges of $\pg(\S)$ incident on $u$. By the discussion preceding Lemma \ref{pinwheel}, each edge of $\st(u)$ is contained in exactly one of $\kk(\S)$ Farey graphs. The following remark offers a characterisation of when two edges of $\st(u)$ are contained in the same Farey graph, and makes apparent that such property is preserved by locally injective simplicial maps.  The proof is immediate. 

\begin{lemma}
Let $u$ be a vertex of $\pg(\S)$. Two edges $e, e' \in \st(u)$ are contained in the same Farey graph if and only if there exists a sequence of edges $e=e_0, e_1, \ldots, e_n=e'$ in $\st(u)$ such that $e_i$ and $e_{i+1}$ are
edges of the same circuit of length 3 in $\pg(\S)$, for all $i = 0, \ldots, n-1$. $\square$ 
\label{labels}
\end{lemma}

We end this subsection with the following observation, which asserts that if a Farey graph $F$ intersects  a stratum  in the pants graph, then either $F$ is contained in the stratum  or else is ``transversal" to it.

\begin{lemma}
Let $F$ be a Farey graph in $\pg(\S)$ and let $T$ be a multicurve. Suppose that $F \cap \pg_T$ has at least 2 vertices. Then $F \subseteq \pg_T$. 
\label{farey-strata}
\end{lemma}

\noindent{\bf{Proof.}} Lemma \ref{farey} implies that $F = \pg_Q$, for some deficiency 1 multicurve $Q$. 
Suppose there exist two distinct vertices $u,v$ in $\pg_Q \cap \pg_T$. Write $u = Q \cup \alpha$, $v= Q \cup \beta$, noting $\alpha \neq \beta$. Since $u,v \in \pg_T$ then $T \subseteq Q$, and therefore  $\pg_Q \subseteq \pg_T$. $\square$


\subsection{Admissible tuples in the pants graph.}

We now introduce the notion of  {\em{admissible tuple}} in the pants graph, a slight generalisation of what Margalit refers to as ``alternating circuit" in \cite{margalit}. 
\begin{definition}[Admissible tuple] Let $n> 3$ and let  $(v_1, \ldots, v_n)$ be a cyclically ordered $n$-tuple of distinct vertices of $\pg(\S)$. We say that $(v_1, \ldots, v_n)$ is {\em{admissible}} if
$v_i$ and $v_{i+1}$ belong to the same Farey graph $F_i$,  and $F_i \neq F_{i+1}$ (counting subindices modulo $n$).	
\end{definition}	

In particular $F_i \cap F_{i+1} = \{v_{i+1}\} $. An admissible 5-tuple will be called a {\em{pentagon}} if $v_i$ and $v_{i+1}$ are adjacent for all $i$. The following lemma describes the structure of admissible $4$- and $5$-tuples in the pants graph, and will be crucial in the proof of our main results. We remark that  this result is implicit in \cite{margalit}.

\begin{figure}[h]
\leavevmode \SetLabels
\L (.10*-.08) $  v_1 $\\
\L (.105*-.17) $ \shortparallel$\\
\L (.04*-.27) $ (\alpha_1, \alpha_2) \cup T  $\\
\L (.04*.94) $(\alpha_4,\alpha_1) \cup T  $\\
\L (.105*.84) $\shortparallel $\\
\L (.10*.75) $v_4  $\\
\L (.27*-.27) $(\alpha_2,\alpha_3) \cup T $\\
\L (.335*-.17) $\shortparallel $\\
\L (.33*-.08) $v_2$\\
\L (.27*.94) $(\alpha_3,\alpha_4) \cup T $\\
\L (.335*.84) $\shortparallel $\\
\L (.33*.75) $v_3$\\
\L (.71*1.04) $ v_4  $\\
\L (.74*1.04) $=$\\
\L (.77*1.04) $ (\alpha_4,\alpha_5)\cup T$\\
\L (.56*-.08) $ v_1 $\\
\L (.565*-.17) $ \shortparallel $\\
\L (.50*-.27) $ (\alpha_1, \alpha_2) \cup T$\\
\L (.87*-.08) $v_2 $\\
\L (.875*-.17) $ \shortparallel $\\
\L (.81*-.27) $ (\alpha_2, \alpha_3) \cup T $\\
\L (.53*.65) $ v_5 $\\
\L (.535*.55) $ \shortparallel $\\
\L (.47*.45) $ (\alpha_5, \alpha_1)\cup T $\\
\L (.90*.65) $v_3$\\
\L (.905*.55) $\shortparallel $\\
\L (.84*.45) $ (\alpha_3,\alpha_4)\cup T $\\
\endSetLabels
\par
\begin{center}
\AffixLabels{\centerline{\epsfig{file =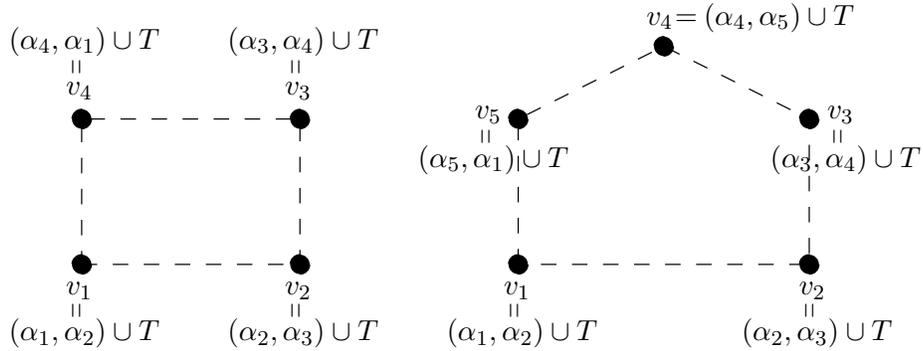, width=10cm, angle= 0}}}
\end{center}

\caption{Admissible 4- and 5-tuples in the pants graph. Here $T$ is a deficiency 2 multicurve, and the dashed line between $v_i$ and $v_{i+1}$ represents a path between
them, entirely contained in a Farey graph $F_i$.}
\label{fig:tuples}
\end{figure}

 \begin{lemma}[Structure of admissible $4$- and $5$-tuples]
 Let  $(v_1, \ldots, v_n)$ be an admissible $n$-tuple, where $n\in\{4,5\}$. Then there exists a deficiency 2 multicurve $T$ such that $v_i \in \pg_T$ for all $i$. Moreover, if $n=4$ then $\S - T$ has exactly 2 non-trivial components, each of complexity 1; if $n=5$ then $\S-T$ has exactly 1 non-trivial component, which has complexity 2. 
 
 \label{structure}
  \end{lemma}

 \noindent{\bf{Proof.}} For the first part, note there is nothing to show if $\kk(\S)=2$, for in that case we let $T = \emptyset$, so that $\pg_T = \pg(\S)$.  So assume $\kk(\S) \geq 3$. Since $v_1, v_2, v_3$ do not belong to the same Farey graph, then  $T= v_1 \cap v_2 \cap v_3$ is a deficiency 2 multicurve. Now $T \subset v_4$ as well; otherwise $v_1$ and $v_4$ would differ by 3 curves and  thus one could not connect $v_1$ and $v_4$ by a path entirely contained in at most 2 Farey graphs.  Similarly, $T \subset v_5$ in the case of an admissible $5$-tuple, and so the first part of the result follows.

Note that, in particular, one can write $v_i=(\aa_i, \aa_{i+1}) \cup T$ for all $i$, as in Figure \ref{fig:tuples}.
Since $T$ has deficiency 2, then $\S - T$ either has one non-trivial component of complexity 2, or  two non-trivial components of complexity 1. Let $X$ (resp. $Y$) be the complexity 1 subsurface filled by $\aa_1$ and $\aa_3$ (resp. $\aa_2$ and $\aa_4$), noting $X \neq Y$ and $X,Y \subset \S -T$.

If $n=4$ then $i(\aa_j, \aa_2) = i(\aa_j, \aa_4) = 0$ for $j \in \{1, 3\}$ (see Figure \ref{fig:tuples}). Thus $X$ and $Y$ are disjoint and thus the result follows. 

Now assume $n=5$. We claim that $X \cup Y$ is connected. If not, then $\aa_5$ is contained in either $X$ or $Y$. But $\aa_5$ is distinct and disjoint from both $\aa_1$ and $\aa_4$ (see Figure \ref{fig:tuples}), and so $\aa_5$ cannot be contained in either $X$ or $Y$, which is a contradiction.  $\square$

\bigskip

 Let $F_1, F_2$ be distinct Farey graphs of $\pg (\S)$ intersecting at a vertex $u$, where $F_i$ is determined by the deficiency 1 multicurve $Q_i$, for $i=1,2$. If the non-trivial components of $\S- Q_1$ and $\S - Q_2$ are disjoint, we say that $F_1$ and $F_2$ {\em{commute}}. Now, Lemma \ref{structure} implies that if $(v_1, \ldots, v_4)$ is an admissible 4-tuple, then $F_i$ and $F_{i+1}$ commute for all $i$, where $F_i$ is the Farey graph containing $v_i$ and $v_{i+1}$. The following converse is an immediate consequence of the proof of Lemma \ref{structure}:

 \begin{corollary}
 Let $F_1, F_2$ be  Farey graphs that commute. For any $u_i \in F_i - \{u\}$, the vertices $ u, u_1, u_2$ are elements of an admissible $4$-tuple.
 \label{commute}
 \end{corollary}
 
 \noindent{\bf{Proof.} }Let $u_i \in F_i-\{u\}$ for $i=1,2$. Since $F_1 \neq F_2$ then $u = (\aa_1, \aa_2) \cup T$, $u_1 = (\aa'_1, \aa_2) \cup T$
 and $u_2 = (\aa_1, \aa_2') \cup T$, for some deficiency 2 multicurve $T$. Now $F_1$ and $F_2$ commute, and so $i(\aa_1', \aa_2') = 0$. Therefore, $(u_2,u,u_1,w)$ is an admissible 4-tuple, where
$w = (\aa_1', \aa_2') \cup T$. $\square$

\bigskip

 In particular, there exists an admissible $4$-tuple in $\pg(\S)$ if and only if $\kk(\S) \geq 3$.  The next technical result will be very important in the next section:

\begin{lemma}[Extending adjacent vertices to admissible tuples]
Let $\S$ be a surface of complexity at least 2. Let $u$ and $v$ be adjacent vertices of $\pg(\S)$ and let $G$ be a Farey graph containing $v$ but not $u$. Then there exists $n \in \{4,5\}$ and a vertex $w \in G -\{v\}$, such that $u, v, w$ are elements of an admissible $n$-tuple.
\label{extension}
\end{lemma}

\noindent{\bf{Proof.}} Write $k = \kk(\S)$, $u = (\aa_1, \aa_2, \dots, \aa_k)$ and $v = (\aa'_1, \aa_2, \dots, \aa_k)$. 
Let $F$ be the Farey graph containing $u$ and $v$, and hence determined by $u \cap v= v - \aa_1'$. Since $u \notin G$ then, up to relabeling the curves of $v$, $G$ is determined by $v - \aa_2$. Let $T$ denote the deficiency 2 multicurve $v - (\aa_1', \aa_2)$. There are two cases to consider:

{\bf{Case 1.}} {\em{$F$ and $G$ commute.}} In this case the result follows from Corollary \ref{commute} by considering any $w \in G$, with $w \neq v$.

{\bf{Case 2.}}  {\em{$F$ and $G$ do not commute.}} Then $\S - T$ has exactly one non-trivial component $S$, of complexity 2, and $(\aa_1, \aa_2)$ and $(\aa_1', \aa_2)$ are adjacent vertices of $\pg(S)$. There are two possibilities for $S$, namely $S$ is a 5-holed sphere or $S$ is a 2-holed torus. 

\begin{figure}[h]
\leavevmode \SetLabels
\L (.26*.3) $\aa $\\
\L (.23*.6) $\bb $\\
\L (.31*.62) $\eta $\\             
\L (.21*.42) $\delta $\\
\L (.31*.42) $\gamma $\\
\L (.505*.28) $\aa $\\
\L (.73*.00) $\bb $\\
\L (.56*-.0) $\gg $\\             
\L (.655*.08) $\delta $\\
\L (.895*.51) $\eta $\\
\endSetLabels
\par
\begin{center}
\AffixLabels{\centerline{\epsfig{file =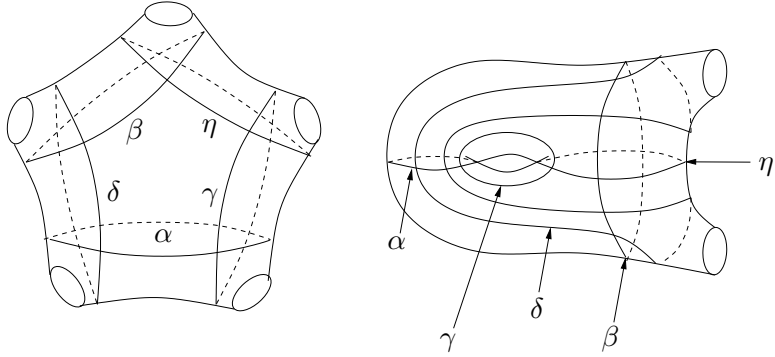, width=10cm, angle= 0}}}
\end{center}
\caption{Curves giving rise to an admissible 5-tuple $(w_1, \ldots, w_5)$ in a 5-holed sphere (left) and a 2-holed torus (right), where $w_1 = (\aa, \bb)$, $w_2 = (\bb,\gg)$, $w_3 = (\gg, \delta)$, $w_4 = (\delta, \eta)$ and $w_5 = (\eta, \aa)$.}
\label{5hs}
\end{figure}

If $S$ is a 5-holed sphere then, up to the action of $\Mod(S)$, the curves $\aa_1, \aa_2, \aa_1'$ are, respectively, the curves $\aa, \bb, \gg$ on the left of Figure \ref{5hs}. Consider the curves $\delta$ and $\eta$, also from the left of Figure \ref{5hs}, and set $v_i = w_i \cup T$, where $w_i$ is defined as in Figure \ref{5hs}. Then $(v_1, \ldots, v_5)$ is an admissible 5-tuple (in this case, a pentagon) in $\pg(\S)$, noting  $v_1 = u$, $v_2 = v$ and $v_3 \in G$. Thus we can take $w = v_3$ and so the result follows. 

The case of $S$ a 2-holed torus is dealt with along the exact same lines, using the curves on the 2-holed torus of Figure \ref{5hs}; in this case, the 5-tuple we obtain is not longer a pentagon (in fact, there are no pentagons in the pants graph of the 2-holed torus; see the proof of Lemma 8 in \cite{margalit})  $\square$

\section{Proof of Theorem B}
\label{proof}

The main ingredient in the proof of Theorem B will be that (locally) injective simplicial maps of pants graphs preserve Farey graphs  and admissible tuples. Let us briefly comment on this. First, a quick argument involving the dual tree of a Farey graph shows that a (locally) injective simplicial map of a Farey graph to itself is in fact bijective; this argument is recorded in Lemma 15 of \cite{shackleton}, and we do not include it here. 

Let $\pp: \pg(\S_1) \to \pg(\S_2)$ be a (locally) injective simplicial map. By the discussion above, if $F$ is a Farey graph, then so is $\pp(F)$. Now, Lemma \ref{labels} implies that if $F, F'$ are distinct Farey graphs that intersect at a vertex, then the same is true for the Farey graphs $\pp(F)$ and $\pp(F')$. 
Using this, plus the definition of admissible tuple, we get that $\pp$ maps admissible $n$-tuples to admissible $n$-tuples.

We can now prove Theorem B. Recall the statement:

\begin{thmb}
Let $\S_1$ and $\S_2$ be compact orientable surfaces such that every connected component of $\S_1$  has positive complexity.  Let $\pp: \pg(\S_1) \rightarrow \pg(\S_2)$ be an injective simplicial map. Then the following hold:

\begin{enumerate}

\item[(1)] $\kk(\S_1) \leq \kk(\S_2)$,

\item[(2)] There exists a multicurve $Q$ on $\S_2$, of cardinality $\kk(\S_2) - \kk(\S_1)$, such that $\pp(\pg(\S_1))= \pg_Q$. In particular, $\pg(\S_1) \cong \pg(\S_2 -Q)$;

\item[(3)] $\S_1$ and $\S_2 -Q$ have the same number of non-trivial components. Moreover, if $X_1, \ldots, X_r$ and $Y_1, \ldots, Y_r$ are, respectively, the non-trivial components of $\S_1$ and $\S_2 - Q$ then, up to reordering the indices, $\pp$ induces an isomorphism $\pp_i: \pg(X_i) \rightarrow \pg(Y_i)$. In particular, $\kk(X_i) = \kk(Y_i)$.

\end{enumerate}

\end{thmb}

\bigskip

\noindent{\bf{Proof.}} Observe that if $\kk(\S_1)=1$ then the result follows from Lemma \ref{farey}. Therefore, from now on we will assume that $\kk(\S_1) \geq 2$.  Let $\kk_i = \kk(\S_i)$, for $i=1,2$.

Part (1) is immediate, since we know $\pp$ maps distinct Farey graphs containing a vertex $u$ (there are $\kk_1$ of these, by Lemma \ref{pinwheel}) to distinct Farey graphs containing $\phi(u)$ (there are $\kk_2$ of these). We will now prove part (2). For clarity, its proof will be broken down into 3 separate claims.

\medskip

\noindent{\bf{Claim I.}} {\em{ Let $u$ be a vertex of $\pg(\S_1)$. There exists a multicurve $Q(u)$ on $\S_2$, of cardinality $\kk_2 - \kk_1$, such that $\pp(e) \subset \pg_{Q(u)}$ for all $e \in \st(u)$.}}

\smallskip

\noindent {\em{Proof.}} By Lemma \ref{pinwheel}, there are  $\kk_1$ Farey graphs  $F_1, \ldots, F_{\kk_1}$ containing $u$. By Lemma \ref{farey}, $\pp(F_i) = \pg_{Q_i}$, for some deficiency 1 multicurve $Q_i \subset \pp(u)$. Consider the multicurve $Q(u) = Q_1 \cap \cdots \cap Q_{\kk_1}$, which has  cardinality $\kk_2 - \kk_1$. Since $Q(u) \subset Q_i$ then  $\pp(F_i) = \pg_{Q_i} \subset \pg_{Q(u)}$. In particular, $\pp(e) \subset  \pg_{Q(u)}$ for all $e \in \st(u)$, since $e \subset F_i$ for some $i$. $\diamond$

\medskip

\noindent{\bf{Claim II.}} {\em{If $v$ is adjacent to $u$, then $\pp(e') \subset \pg_{Q(u)}$ for all $e' \in \st(v)$, where $Q(u)$ is the multicurve given by Claim I for $u$.}}

\smallskip

\noindent{\em{Proof.}} Let $e$ be the edge with endpoints $u$ and $v$. Let $e' \in \st(v)$ and let $F$ be the unique Farey graph containing $e'$. If $e \subset F$ then the result follows from the proof of Claim I. So suppose $e$ is not contained in $F$, so $u \notin F$. By Lemma \ref{extension}, there exist a number $n\in\{4,5\}$ and a vertex $w \in F$, with $w \neq v$, such that $u,v,w$ are elements of an admissible $n$-tuple in $\pg(\S_1)$,  which we denote by $\tau$. Thus  $\pp(u), \pp(v), \pp(w)$ are also elements of an admissible $n$-tuple in $\pg(\S_2)$. Therefore there is a deficiency 2 multicurve $T$ on $\S_2$ such that $\pp(\tau) \subset \pg_T$, using Lemma \ref{structure}.

We now claim that $Q(u) \subseteq T$. To see this, let $z$ be the unique element of $\tau - \{v\}$ which is contained in the same Farey graph as $u$, noting $z,u,v$ are not contained in the same Farey graph of $\pg(\S_1)$ by the definition of admissible tuple. Therefore  $\phi(z), \phi(u),\phi(v)$ are not contained in the same Farey graph of $\pg(\S_2)$ and so $\phi(z) \cap \phi(u) \cap \phi(v) = T$, since $\phi(\tau) \subset \pg_T$ and $T$ has deficiency 2. Finally, $Q(u) \subseteq \phi(z) \cap \phi(u) \cap \phi(v)$ since $\phi$ maps every Farey graph containing $u$ into $\pg_{Q(u)}$ and $u, v$ (resp. $u, z$) are contained in a common Farey graph. Thus $Q(u) \subseteq T$, as desired.

Since $Q(u) \subseteq T$ then $\phi(\tau) \subset \pg_T \subseteq \pg_{Q(u)}$. In particular, $\phi(w)$ is contained in $\pg_{Q(u)}
$ and thus in $\phi(F) \cap \pg_{Q(u)}$. Since $\phi(v) \in \phi(F) \cap \pg_{Q(u)}$ as well, we conclude that $\phi(F) \subset \pg_{Q(u)}$ by Lemma \ref{strata}. In particular, $\pp(e') \subset \pg_{Q(u)}$ and thus Claim II follows. $\diamond$

\medskip

As a consequence, and since $\pg(\S_1)$ is connected, it follows that  $\pp(\pg(\S_1)) \subseteq \pg_Q$, where $Q = Q(u)$ for some, and hence any, vertex $u$ of $\pg(\S_1)$. Actually, more is true:

\medskip

\noindent{\bf{Claim III.}} {\em{$\phi(\pg(\S_1)) = \pg_Q$.}}

\smallskip

\noindent{\em{Proof.}} Let $e$ be an edge of $\pg_Q$; we want to show that $e \in {\rm{Im}}(\phi)$. Since $\phi(\pg(\S_1))$ and $\pg_Q$ are connected, and  since $\phi(\pg(\S_1)) \subseteq \pg_Q$, we can assume $e \in \st(\pp(u))$ for some vertex $u$ of $\pg(\S_1)$. Note that $e$ is contained in a unique Farey graph $H$ and that $H \subset  \pg_Q$ by Lemma \ref{strata}.

Since $\pg_Q \cong \pg(\S_2 -Q)$ and $\kappa_1 =  \kappa(\S_2 - Q)$,  there are exactly $\kk_1$ distinct Farey graphs in $\pg(\S_2)$ which are contained in $\pg_Q$ and contain $\phi(u)$, by Lemma \ref{pinwheel}. Again by Lemma \ref{pinwheel}, there are exactly $\kk_1$ distinct Farey graphs in $\pg(\S_1)$ containing $u$. Since $\phi$ maps distinct Farey graphs containing $u$ to distinct Farey graphs containing $\phi(u)$, we get that $H = \phi(F)$ for some Farey graph $F$ in $\pg(\S_1)$ containing $u$. In particular, $e \in {\rm{Im}}(\phi)$, as desired. $\diamond$

\bigskip

Finally, we will prove part (3). Let $X_1, \ldots, X_r$ be the non-trivial components of $\S_1$. Observe that every pants decomposition of $\S_1$ has the form $(v_1, \ldots, v_r)$, where $v_i$ is a pants decomposition of $X_i$, and so  $\pg(\S_1) = \Pi_{i=1}^r\pg(X_i)$.  Fix a pants decomposition $v=  (v_1, \ldots, v_r)$ of $\S_1$. Then $\pp$ induces an injective simplicial map 
\[\pp_{i}: \pg(X_i) \to \pg_Q \cong \pg(\S_2 - Q),\]

\noindent by  setting $\pp_i(w) = \pp(v_1, \ldots, v_{i-1}, w, v_{i+1}, \ldots, v_r)$ for all vertices $w$ of $X_i$. Applying part (2) of Theorem B to $\pp_i$, we deduce that there exists an essential subsurface $Y_i$ of $\Sigma_2 - Q$ such that $\pp_i(\pg(X_i) ) =  \pg(Y_i)$. In particular, $\kk(Y_i) = \kk(X_i)$, by part (1). 
Moreover, by discarding those connected components of $Y_i$ homeomorphic to a 3-holed sphere, we can assume that $Y_i$ has no trivial components. 

\medskip

\noindent{\bf{Claim.}} $Y_i$ is connected. 

\smallskip

\noindent{\em{Proof.}} Suppose, for contradiction, that $Y_i$ had $N\geq 2$ components $Z_1, \ldots, Z_N$. In particular $0<\kk(Z_j) < \kk(Y_i) = \kk(X_i)$, for all $j$, and 
$\pp_i(\pg(X_i) ) =  \pg(Z_1) \times \ldots \times  \pg(Z_N)$. Thus, the image of an edge of $\pg(X_i)$ under $\pp_i$ is contained in one of the factors above, and thus the same holds for the image of any Farey graph under $\pp_i$, by Corollary \ref{commute}. Moreover, if $F$ and $F'$ do not commute, then $\pp_i(F)$ and $\pp_i(F')$ are contained in the same factor, also by Corollary \ref{commute}.

Let $u$ be a vertex of $\pg(X_i)$. We now define the {\em{adjacency graph} $\Gamma$} of $u$,  introduced independently by Behrstock-Margalit  \cite{behrstock-margalit} and Shackleton \cite{shackleton}. The vertices of $\Gamma$ are exactly those curves in $u$, and two distinct curves are adjacent in $\Gamma$ if they are boundary components of the same pair of pants determined by $u$. Observe that $\Gamma$ is connected since $X_i$ is. 

Now a Farey graph containing $u$ is  is determined by a deficiency 1 multicurve contained in $u$ or, equivalently, by a curve in $u$.  Moreover, two curves in $u$ are adjacent if and only if the Farey graphs they determine do not commute.  Let $\cal{F}$ be the graph whose vertices are those Farey graphs containing $u$ and whose edges correspond to distinct non-commuting Farey graphs.  Note $\cal{F}$ is isomorphic to $\Gamma$ and so it is connected. 

By Lemma \ref{pinwheel} and since $\cal{F}$ is connected, there exist $\kk(X_i)$ Farey graphs in $\pg(X_i)$, all containing $u,$ which are mapped into the same factor of $\pg(Z_1) \times \ldots \times  \pg(Z_N)$ under $\pp_i$. This contradicts Lemma \ref{pinwheel}, since $\kk(Z_j) < \kk(X_i)$ for all $j$, and thus the claim follows.  $\diamond$

\bigskip

The discussion above implies that there are $r$ connected subsurfaces $Y_1, \ldots,Y_r$ of $\S_2 - Q$ such that, up to reordering, $\phi$ induces an isomorphism $\pp_i: \pg(X_i) \to \pg(Y_i)$ for $i = 1, \ldots, r$. In particular, $\kk(X_i) = \kk(Y_i)$. Now,

\[\S_1 = X_1 \sqcup \ldots \sqcup X_r\] and

\[\S_2 - Q \supseteq Y_1 \cup \ldots \cup Y_r,\]

\noindent and therefore the $Y_i$'s are pairwise disjoint, since the $X_i$'s are pairwise disjoint, $\kk(X_i) = \kk(Y_i)$ and $\kk(\S_1) = \kk(\S_2 - Q) $. For the same reason, they are the only non-trivial connected components of $\S_2 - Q$.
This finishes the proof of Part (3) of Theorem B. $\square$

\vspace{1cm}

\section{A consequence of Theorem B}
\label{apps}

We now present an application of Theorem B to pants graph automorphisms. Let $\phi: \pg(\S) \to \pg(\S)$ be an injective simplicial map; by Theorem B, $\phi$ is in fact an isomorphism.  Let $\aa$ be a curve on $\S$, and observe that $\pg(\S-\aa) \cong \pg_\aa \subset \pg(\S)$. Then $\phi$ induces an injective simplicial map, which we also denote by $\phi$, from $\pg_\aa$ to $\pg(\S)$. Applying Theorem B to $\S_1 = \S - \alpha$ and  $\Sigma_2 = \Sigma$, we readily obtain the following corollary, which implies that pants graph automorphisms preserve the pants graph stratification:

\begin{corollary}
Let $\pp:\pg(\S) \rightarrow \pg(\S)$ be an automorphism. Then, for every curve $\aa$, there exists a unique curve $\bb$ such that $\pp(\pg_\aa)=\pg_\bb$. Moreover, $\S- \aa$ and $\S - \bb$ have the same number of non-trivial components. \label{strata}
\end{corollary}

In \cite{margalit}, Margalit introduced the notion of a {\em{marked}} Farey graph in the pants graph. As Farey graphs, they are preserved by pants graph automorphisms. A marked Farey graph singles out exactly one curve on $\S$, although there are infinitely many marked Farey graphs in $\pg(\S)$ that single out a given curve. Margalit associates, to the pants graph automorphism $\pp$, a curve graph automorphism $\psi$ by defining $\psi(\aa)$ to be the curve $\bb$ singled out by $\pp(F)$, where $F$ is a marked Farey graph that singles out $\aa$. One of the main steps in \cite{margalit} is to show that this construction gives rise to a well-defined map between the pants graph automorphism group and the curve graph automorphism group, which Margalit then shows is an isomorphism.  If $\S$ is not a 2-holed torus, Theorem \ref{dan} then follows from results of Ivanov \cite{ivanov}, Korkmaz \cite{korkmaz} and Luo \cite{luo} on the automorphism group of the curve graph. The case of the 2-holed torus requires separate treatment in \cite{margalit}, and boils down to showing that the curve graph automorphism induced by a pants graph automorphism maps (non-)separating curves to (non-)separating curves.

Similarly, one could define a curve graph automorphism $\psi$ from the pants graph automorphism $\pp$, by setting $\psi(\aa) = \bb$, where $\bb$ is the curve such that $\pp(\pg_\aa) = \pg_\bb$ in Corollary \ref{strata}. One quickly checks that this produces an isomorphism between the pants graph automorphism group and the curve graph automorphism group, and thus Theorem \ref{dan} follows if the surface is not the 2-holed torus. The case of the 2-holed torus is also deduced from Corollary \ref{strata} by applying the exact same argument as in \cite{margalit}. We remark that this approach to pants graph automorphisms is similar in spirit to those of Masur-Wolf \cite{masur-wolf} and Brock-Margalit \cite{brock-margalit} for showing that Weil-Petersson isometries are induced by surface self-homeomorphisms. Indeed, one of the key steps there is to prove that Weil-Petersson isometries preserve the stratification of the Weil-Petersson completion.


\section{Proof of Theorem A}
\label{final}

We are finally ready to give a proof of Theorem A. We will need the following lemma, which we believe is folklore, but which we nevertheless prove for completeness. 

\begin{lemma}[Classification of pants graphs up to isomorphism]
Let $\S$, $\S'$ be two compact connected orientable surfaces of complexity at least 2. Then $\pg(\S)$ and $\pg(\S')$ are isomorphic if and  only if $\S$ and $\S'$ are homeomorphic. 
\label{isom-class}
\end{lemma}

\noindent{\bf{Proof.}} First, by Part (1) of Theorem B, if $\pg(\S)$ and $\pg(\S')$ are isomorphic then $\kk(\S) = \kk(\S')$. We consider the following three cases:

\medskip

 (i) Suppose $\kk(\S) > 3$. By Theorem \ref{dan}, ${\rm{Aut}}(\pg(\S)) \cong {\rm{Mod}}(\S)$. Thus if 
 $\pg(\S) \cong \pg(\S')$ then ${\rm{Mod}}(\S) \cong {\rm{Mod}}(\S')$ and thus $\S$ and $\S'$ are homeomorphic by Theorem 2 of \cite{shackleton}. (We remark that Shackleton's result was first proved by Ivanov-McCarthy \cite{ivanov-mccarthy} for surfaces of positive genus.)

\medskip

(ii) Suppose that $\kk(\S) = 2$. Up to renaming the surfaces, $\S$ is a 5-holed sphere and $\S'$ is a 2-holed torus. The curves on the 5-holed sphere on the left of Figure \ref{5hs} yield
the existence of a pentagon in $\pg(\S)$, while there are no pentagons in the pants graph of the 2-holed torus (see the proof of Lemma 8 in \cite{margalit}).

\medskip

(iii) Finally,  we consider the case $\kk(\S)=3$. Let us denote by $S_{g,b}$ the surface of genus $g$ with  $b$ boundary components. We have that 

\[(\S, \S') \in \{ (S_{0,6}, S_{2,0}), (S_{1,3}, S_{0,6}), (S_{1,3}, S_{2,0})\} \]

Suppose, for contradiction, that there exists an isomorphism $\pp$ between $\pg(S_{0,6})$ and $\pg(S_{2,0})$. Choose a curve $\aa$ on $S_{0,6}$ such that $S_{0,6} - \aa = S_{0,3} \sqcup S_{0,5}$, noting that $\pg_\aa \cong \pg(S_{0,5})$. By Theorem B, there exists a curve $\bb$ on $S_{2,0}$ such that $\pp(\pg_\aa) = \pg_\bb$. Moreover, $S_{2,0} - \bb$ has to be connected, and thus $S_{2,0} - \bb$ is homeomorphic to $S_{1,2}$. In particular, $\pg(S_{2,0} - \bb) \cong \pg(S_{1,2})$. Thus we get an isomorphism between $\pg(S_{0,5})$ and $\pg(S_{1,2})$, which contradicts (ii). 

Now, suppose there is an isomorphism $\pp$ between $\pg(S_{1,3})$ and $\pg(S_{0,6})$. We choose a curve $\aa$ on $S_{1,3}$ such that $S_{1,3} - \aa = S_{0,3} \sqcup S_{1,2}$. Arguing as above, 
we get an isomorphism between $\pg(S_{1,2})$ and  $\pg(S_{0,5})$, which does not exist, by (ii).

Finally, suppose that there exists an isomorphism  between $\pg(S_{1,3})$ and $\pg(S_{2,0})$. Choose
a non-separating curve $\gg$ on $S_{1,3}$, so that $S_{1,3} - \gg = S_{0,5}$. By the same argument as above, we get an isomorphism between $\pg(S_{0,5})$ and  $\pg(S_{1,2})$, contradicting (ii). $\square$.

\bigskip

As  mentioned before, the pants graph of the two complexity 1 surfaces (that is, the 1-holed torus and the 4-holed sphere) are isomorphic. We remark that the isomorphism classification of pants graphs is slightly different to that of curve complexes (see Lemma 2.1 in \cite{luo}).

\bigskip

\noindent{\bf{Proof of Theorem A.}} Let $\pp: \pg(\S_1) \rightarrow \pg(\S_2)$ be an injective simplicial map. By Theorem B, there exists a multicurve $Q$ on $\S_2$, of deficiency $\kk(\S_1)$, such that $\pp(\pg(\S_1)) =  \pg_Q \cong \pg(\S_2 - Q)$. Discarding the trivial components of $\S_2 - Q$ we obtain an essential subsurface $Y \subset \S_2 - Q$, with no trivial components, and such that  $\pg(\S_1) \cong \pg(Y)$.
We can thus view $\phi$ as an isomorphism  $\pg(\S_1) \to \pg(Y)$. Let us first assume that $\S_1$ is connected. In that case $Y$ is connected as well, by part (3) of Theorem B. Since $\kk(\S_1) \geq 2$, Lemma \ref{isom-class} implies there exists a homeomorphism $g: \S_1 \to Y$, which induces an isomorphism $\psi: \pg(\S_1) \to \pg(Y)$ by $\psi(v) = g(v)$. By Theorem \ref{dan}, there exists $f \in \Mod(Y)$ such that $\phi = f \circ \psi$. Thus $f \circ g$ induces $\phi$. 

If $\S_1$ is not connected, the result follows by applying the above argument to each connected component of $\S_1$. $\square$


\bigskip

\begin{small}
\noindent Javier Aramayona\\
National University of Ireland, Galway\\
\url{Javier.Aramayona@nuigalway.ie}
\end{small}

\end{document}